\documentstyle[12pt]{article}
\newcommand{\braid}[2]{{#1}$\lower4pt\hbox{$\oo\atop\raise4pt\hbox{$\scripts

criptstyle\Phi $}$}${#2}}
\newcommand{\twist}[2]{{#1}${\,\scriptscriptstyle
\Phi}\atop\raise9pt\hbox{$\scriptstyle\oo$} ${#2}}

\newcommand{\ve}{\varepsilon}
\newcommand{\z}{{\zeta}}

\newcommand{\barr}{\begin{displaymath} \begin{array}{lll}}
\newcommand{\earr}{\end{array} \end{displaymath}}
\newcommand{\be}{\begin{eqnarray}}
\newcommand{\ee}{\end{eqnarray}}
\newcommand{\n}{\nonumber }
\newcommand{\oo}{\otimes}
\newcommand{\al}{\alpha}
\newcommand{\bt}{\beta}

\newcommand{\si}{\sigma}
\newcommand{\gm}{\gamma}
\newcommand{\h}{\eta}

\newcommand{\lmx}[1]{\begin{displaymath} {#1}= \left(\begin{array}{rrrr}}
\newcommand{\rmx}{\end{array} \right), \end{displaymath}}
\newcommand{\rmtx}{\end{array} \right). \end{displaymath}}

\newcommand{\lmxx}[1]{\begin{displaymath} {#1}\rightarrow
\left(\begin{array}{rr}}
\newcommand{\rmxx}{\end{array} \right), \end{displaymath}}

\newtheorem{theor}{Theorem}

\newtheorem{lemma}{Lemma}

\def\C{{\mathchoice {\setbox0=\hbox{$\displaystyle\rm C$}\hbox{\hbox
to0pt{\kern0.4\wd0\vrule height0.9\ht0\hss}\box0}}
{\setbox0=\hbox{$\textstyle\rm C$}\hbox{\hbox
to0pt{\kern0.4\wd0\vrule height0.9\ht0\hss}\box0}}
{\setbox0=\hbox{$\scriptstyle\rm C$}\hbox{\hbox
to0pt{\kern0.4\wd0\vrule height0.9\ht0\hss}\box0}}
{\setbox0=\hbox{$\scriptscriptstyle\rm C$}\hbox{\hbox
to0pt{\kern0.4\wd0\vrule height0.9\ht0\hss}\box0}}}}
\def\R{{\mathchoice {\setbox0=\hbox{$\displaystyle\rm R$}\hbox{\hbox
to0pt{\kern0.4\wd0\vrule height0.9\ht0\hss}\box0}}
{\setbox0=\hbox{$\textstyle\rm R$}\hbox{\hbox
to0pt{\kern0.4\wd0\vrule height0.9\ht0\hss}\box0}}
{\setbox0=\hbox{$\scriptstyle\rm R$}\hbox{\hbox
to0pt{\kern0.4\wd0\vrule height0.9\ht0\hss}\box0}}
{\setbox0=\hbox{$\scriptscriptstyle\rm R$}\hbox{\hbox
to0pt{\kern0.4\wd0\vrule height0.9\ht0\hss}\box0}}}}

\begin{document}
\begin{titlepage}
\begin{center}
{\Large \bf Inhomogeneous quantum Lie algebras}
\end{center}
\vspace{0.5cm}
\begin{center}
 {\bf P. P. Kulish\footnote{Partially supported by the RFFI
 grant 98-01-00310.}}
\end{center}
\vspace{0.5cm}
\begin{center}
St.Petersburg Department of the Steklov
Mathematical Institute,\\
Fontanka 27, St.Petersburg, 191011,
Russia \\ (kulish@pdmi.ras.ru)
\end{center}
\vspace{0.5cm}
\begin{center}
 {\bf A. I. Mudrov}
\end{center}
\vspace{0.5cm}
\begin{center}
Department of Theoretical Physics,
Institute of Physics, St.Petersburg State University, Ulyanovskaya 1,
St.Petergof, St.Petersburg, 198904, Russia\\
(aimudrov@dg2062.spb.edu)
\end{center}
\vspace{1cm}
\begin{center}
{\bf Abstract}
\end{center}
{\footnotesize
We study quantization of a class of inhomogeneous
Lie bialgebras which are crossproducts in
dual sectors with Abelian invariant parts.
This class forms a category stable under dualization and the double
operations. The quantization turns out to be a functor commuting with
them. The Hopf operations and the universal R-matrices are
given in terms of generators. The quantum algebras obtained appear
to be isomorphic to the universal enveloping Poisson-Lie
algebras on the dual groups.
}
\vspace{0.5in}

\begin{center}
{1999}
\end{center}

\end{titlepage}
\section{Introduction}
Inhomogeneous Lie groups such as
those belonging to the Cayley-Klein series, including
Poincar\'e and  Galilei, play important
role in classical physics and geometry. They realize the maximal
sets of (continuous) symmetries of the simply connected
(pseudo) Riemannian spaces of the zero curvature.
The generalization of the semi-direct product of classical
groups in the framework of non-commutative geometry is
the bicrossproduct \cite{Mj4} of two Hopf algebras ${\cal A}$ and
${\cal B}$ characterized by actions of
${\cal A}$ on ${\cal B}$ and  ${\cal B}^*$ on ${\cal A}^*$.
Nowadays, there are numerous examples of bicrossproducts  known,
including those among quantum deformations of the
Cayley-Klein algebras \cite{CK}. Unfortunately, contractions of
the quantum orthogonal algebras leading to those solutions
result in poles in their classical r-matrices, however
disappearing from the skew-symmetric part. Thus a Lie bialgebra
survives, whereas the quasitriangular structure is broken.
The canonical (and the simplest) example
of the bicrossproduct construction is
the second (non-standard) quantization of the Borel subalgebra
$b(2)\subset sl(2)$. At the same time, this algebra is the
result \cite{Ogiv} of Drinfeld's twist \cite{D1,D2}
of the universal enveloping algebra $U(b(2))$.
Another examples of twisted bicrossproduct Hopf algebra are the
null-plane quantized Poincar\'e algebra
\cite{AHO} and extended jordanian deformations
of $U(sl(N))$. These quantizations involve special
non-degenerate 1-cocycles on Lie groups \cite{M12,KLM,ES}.
All those algebras are twist-equivalent to classical universal enveloping
algebras, and that equivalence holds for their representation theories.
Quasitriangular  bicrossproduct Hopf algebras with non-unitary
R-matrices were found in \cite{Mj2}
via the quantum double construction in the framework of the
matched pairs of finite groups. The present work is devoted to the
study of "continuous" bicrossproduct Hopf algebras
with Abelian invariant subalgebras. In the classical
differential geometry these correspond to inhomogeneous
Lie groups, containing sets of commutative translations.
Quantum version of the theory appears to possess a number of
remarkable features,
for example, invariance of the category of interest with
respect to dualization and the double procedures.
Explicitly built Hopf operations  allows us to conduct
the detailed study  of quantum  doubles,
construct canonical elements and R-matrices for generic
quasitriangular Hopf algebras from the category under
investigation.
Our approach relies on a kind of "universality" of
the double construction, proved by Radford \cite{Rad}
and meaning the following.
Every quasitriangular Hopf algebra contains a minimal quasitriangular
Hopf subalgebra which is a quotient of the quantum double of
another Hopf subalgebra.

The Hopf algebras studied in this paper are related to the
matched pairs of continuous groups, that explains
appearance of objects inherent to classical differential geometry,
such as Lie group 1-cocycles. In fact, quantum commutation
relations turn out to be just the Poisson brackets on
the dual Lie group, and the quantum symmetries form
the universal enveloping algebra of the corresponding
Poisson-Lie algebra of functions.

\section{Quantum double and quasitriangularity}
The purpose of this preliminary section is to present, for completeness,
to prove that every quasitriangular Hopf algebra contains
a subalgebra which is a quotient of the quantum double \cite{Rad}.
We start with the following elementary proposition
from the linear algebra.
\begin{lemma}
\label{CE}
Let ${\bf L}$ be a vector space and $r\in {\bf L}\otimes {\bf L}$.
Consider the subspaces ${\bf L}_+=r({\bf L}^*)$,
$r(x)= \langle x\otimes id,r\rangle$ and ${\bf L}_-=r^*({\bf L}^*)$.
Then ${\bf L}_+^*\sim{\bf L}_-$ and the element $r$
coincides with the image of the canonical element under
the induced map
${\bf L}^*_+\oo {\bf L}_+ \to {\bf L}_-
\oo {\bf L}_+ \subset {\bf L}\oo{\bf L}$
identical on the second tensor factor.
\end{lemma}
The first part of the statement follows from the commutative diagrams

\begin{picture}(400,100)
\put(50,70){${\bf L}^*$}
\put(80,40){${\bf L}_+$}
\put(52,10){$0$}

\put(112,70){${\bf L}$}
\put(112,10){$0$}

\put(62,20){\vector(1,1){15}}
\put(62,67){\vector(1,-1){15}}

\put(95,53){\vector(1,1){15}}
\put(95,35){\vector(1,-1){15}}

\put(65,73){\vector(1,0){42}}

\put(80,78){$r$}
\put(55,50){$\tilde r$}
\put(110,50){$i$}

\put(250,70){${\bf L}$}
\put(280,40){${\bf L}^*_+$}
\put(252,10){$0$}

\put(312,70){${\bf L}^*$}
\put(312,10){$0$}

\put(275,33){\vector(-1,-1){15}}
\put(275,54){\vector(-1,1){15}}

\put(310,68){\vector(-1,-1){15}}
\put(310,20){\vector(-1,1){15}}

\put(307,73){\vector(-1,0){42}}

\put(280,78){$r^*$}
\put(255,50){$\tilde r^*$}
\put(310,50){$i^*$}
\end{picture}

\noindent
where the isomorphism ${\bf L}^*_+\to {\bf L}_-$ is given
by the map $\tilde r^*$. Let us prove that the
map $\tilde r^* \oo id $ brings the canonical
element $l\in {\bf L}^*_+ \oo {\bf L}_+$ right to
$r\in {\bf L}_- \oo {\bf L}_+\subset {\bf L} \oo {\bf L}$.
Indeed, for every $x,y\in {\bf L}^*$ we have
$\langle (\tilde r^*\oo id)(l),x\oo y\rangle =
\langle l,r(x)\oo y\rangle = \langle r(x),y\rangle =
\langle r,x\oo y\rangle$. Here we used the
characteristic of the canonical element,
$\langle l,r(x)\oo id\rangle =r(x)$.

Now consider the quasi-classical situation when  ${\bf L}$ is
a Lie bialgebra. Then ${\bf L}_\pm$ and their linear
sum are themselves sub-bialgebras. Moreover,
${\bf L}_+ +{\bf L}_-$ is the minimal quasitriangular
Lie sub-bialgebra, where the classical r-matrix lives in fact.
Since $r^*$ is a coalgebra homomorphism but an algebra
anti-homomorphism, it also can be regarded as a morphism
in the Lie bialgebra category, ${\bf L}^*_+$ being endowed with
the opposite bracket. Let us consider the double ${\cal D}({\bf L}_+)$
built on the linear sum of ${\bf L}_\pm$ and including them as
Lie sub-bialgebras \cite{D3}.
In an evident way the mapping
${\cal D}({\bf L}_+)\to {\bf L}_++{\bf L}_-\subset {\bf L}$
is defined, which is just the identification
embedding on each addends. Its restrictions on ${\bf L}_\pm$
preserves the Lie structures separately. Let us prove the same assertion with
respect
to the commutator $[{\bf L}_+,{\bf L}_-]$.
For arbitrary $x,y,z \in {\bf L}^*$ consider the classical Yang-Baxter
equation
$$\langle[r_{12},r_{13}],x\otimes y\otimes z\rangle +
\langle[r_{13},r_{23}],x\otimes
y\otimes z\rangle+
 \langle[r_{12},r_{23}],x\otimes y\otimes z\rangle=0.$$
Having introduced the notations $x_+=r(x)\in {\bf L}_+$, $x_-=r^*(x)\in
{\bf L}_-$,
where $x\in {\bf L}^*$, rewrite this equality as
$$\langle[y_-,z_-],x\rangle + \langle[x_+,y_+],z\rangle+
\langle[x_+,z_-],y\rangle=0.$$
This, in its turn, is equivalent to
$$[x_+,z_-]= (x_+\triangleright z)_- - (z_-\triangleright x)_+,$$
where $\triangleright = -ad^*|_{\bf L}$ - is the conjugate to the
adjoint representation. As the mapping
$i^*\colon {\bf L}^*\to {\bf L}^*_+\sim {\bf L}_- $
is a homomorphism of ${\bf L}_+$-modules (and the same is the case
with replacement $\pm \to \mp$), the latter expression
can be rewritten in the form
$$[x_+,z_-]= x_+\triangleright z_- - z_-\triangleright x_+,$$
where $\triangleright $ is already considered as  $-ad^*|_{{\bf L}_\pm}$.
But this is exactly the definition of the Lie bracket in ${\cal D}({\bf L}_+)$

Now the similar result is formulated for abstract quasitriangular
Hopf algebras (strictly speaking, finite dimensional). Recall
that a Hopf algebra ${\cal H} $ is
quasitriangular \cite{D3} if there exists an element  ${\cal R}\in {\cal H}^{\oo2}$
(the universal R-matrix) such
that
$$
(\Delta\oo id) ({\cal R})= {\cal R}_{13}{\cal R}_{23}, \quad
(id\oo \Delta) ({\cal R})= {\cal R}_{13}{\cal R}_{12},
$$
$$
{\cal R}\Delta(h)= \Delta'(h){\cal R},
$$
where the prime denotes the opposite coproduct and the subscripts indicate the
way of embedding into the tensor cube.
It follows from here that ${\cal R}$ satisfies the Yang-Baxter equation
$$
{\cal R}_{12}{\cal R}_{13}{\cal R}_{23}=
{\cal R}_{23}{\cal R}_{13}{\cal R}_{12}.
$$
R-matrix defines two algebra and anti-coalgebra homomorphisms
from ${\cal H}^*$ to ${\cal H}$, $\h\to \langle\h\oo id, {\cal R}\rangle$
and $\h\to \langle id\oo \h, {\cal R}^{-1}\rangle=\langle id\oo \h, (S\oo
id)({\cal R})\rangle$,
their images denoted ${\cal H}_+$ and ${\cal H}_-$. Hopf algebra
${\cal H}_-$ is isomorphic to ${\cal H}^*_+$ taken with
the opposite multiplication.  Let us consider the double \cite{D3}
${\cal D}({\cal H}_+)$, which is built on the tensor product of
${\cal H}_+$ and ${\cal H}_{+,op}^*$ embedded there as sub-bialgebras.
The relations between these two factors are encoded in
the Yang-Baxter equation on the canonical element
$I=h^i\oo h_i\in {\cal H}_{+,op}^*\oo{\cal H}_+$.
The map ${\cal D}({\cal H}_+)\to {\cal H}$ defined as identical
on ${\cal H}_+\oo 1$ and the isomorphism
$1\oo{\cal H}_{+,op}^*\to {\cal H}_-$
respects the bialgebra structures when restricted to
these sub-bialgebras. The image of the canonical element under this
mapping is the R-matrix (cf. Lemma \ref{CE}), and the  cross-relations
in the quantum double go over into the quantum Yang-Baxter
equation fulfilled by the R-matrix. Hence the map of concern
is a homomorphism. Its surjective image includes sub-Hopf algebras
${\cal H}_\pm$ and is exactly that subalgebra in ${\cal H}$
where the R-matrix actually lies. Thus we finish the proof.

The subspaces ${\bf L}_+$ and ${\bf L}_-$ glue over the Cartan subalgebra
in a standard (Drinfeld-Jimbo) semisimple Lie bialgebra.
Its quantization belongs to the class  of the factorizable Hopf algebras
introduced in \cite{RS}.
For that type of algebras, the "universality" property of the double
was stated therein. Alternative examples are triangular bialgebras with
skew-symmetric r-matrices, where ${\bf L}_+$ coincides with ${\bf L}_-$.
The simplest case of the double of the triangular quantized $sl(2)$-Borel
subalgebra was studied in \cite{BH,V}.

\section{Inhomogeneous Lie bialgebras and their quantization}
We introduce the bicrossproduct structure on a Lie bialgebra
by means of an involution $\si$ which is assumed to be  an anti-automorphism
of the algebra and an automorphism of coalgebra, demanding
that $\si$- and $(-\si^*)$- invariant subspaces should be Lie subalgebras
in ${\bf  L}$ and ${\bf  L}^*$ and therefore commutative subalgebras.
The exact description, in terms of generators
$H_i \in {\bf H}$ and  $X^\mu \in {\bf V}$,
${\bf L} ={\bf H}\triangleright {\bf V}$,
is as follows.
\begin{equation}
\begin{array}{lcl}
 [H_i, H_k]&=& C^m_{ik} H_m ,                       \\[6pt]
 [H_i,X^\mu] &=& A_{i\nu}^{\>\>\mu} X^\nu,\\[6pt]
 [X^\mu,X^\nu] &=& 0,\\
 \delta(X^\mu)&=&\gm^\mu_{\rho\si}(X^\rho \oo X^\si),\\
 \delta(H_i)&=&\al_{\rho i}^{\>\>k} (X^\rho\oo H_k- H_k \oo X^\rho).
\end{array}
\label{COMMREL}
\end{equation}
The tensors $C^m_{ik}$ and $\gm^\mu_{\rho\si}$ are skew-symmetric and satisfy
the Jacobi identity. Matrices $A_i$ and $\al_\mu$ realize
representations of ${\bf H}$ on ${\bf V}$ and ${\bf V}^*$ on ${\bf H}^*$,
respectively.

To match Lie bialgebra cocycle condition, Lie structures
on ${\bf  L}$ and ${\bf  L}^*$ should be consistent:
\be
A_{i\nu}^{\>\>\mu} \gm^\nu_{\rho\si}-\gm^\mu_{\nu\si}A_{i\rho}^{\>\>\nu}
-\gm^\mu_{\rho\nu}A_{i\si}^{\>\>\nu}
&=&
A_{k\si}^{\>\>\mu} \al_{\rho i}^{\>\>k}-
A_{k\rho}^{\>\>\mu} \al_{\si i}^{\>\>k},
\label{cond1}\\
\al_{\mu m}^{\>\>k} C^m_{ij}-C^k_{im}\al_{\mu j}^{\>\>m}
    -C^k_{mj}\al_{\mu i}^{\>\>m}
&=&
\al_{\nu j}^{\>\>k} A_{i\mu}^{\>\>\nu} - \al_{\nu i}^{\>\>k}
A_{j\mu}^{\>\>\nu}.
\label{cond2}
\ee
Bialgebras of such a type form a category which we
denote ${\cal B}$. Its morphisms are those respecting
Lie products on ${\bf  L}$ and ${\bf  L}^*$ and commuting
with the involution $\si$.

Our quantization strategy relies on
the quantum duality principle \cite{STS9,D3} as applied to the
problem of "exponentiating" bialgebras of concern.
Following this  principle, we consider
a quantum algebra as a variety of noncommutative functions on
the group $\exp({\bf L}^*)$. In accordance with the dual group method of
building quantum deformations \cite{LM}, we fix the coproduct
\begin{equation}
\begin{array}{lcl}
 \Delta (1) &=& 1\oo 1,\\
 \Delta (X^\mu) &=& D^\mu(X\oo 1,1\oo X)\\
 \Delta (H_i) &=& (e^{\al \cdot X})^k_i \oo H_k  + H_i \oo 1.
\end{array}\label{coproduct:1}
\end{equation}
just exponentiating the Lie bracket on ${\bf L}^*$.
We use notations $D(.,.)$ for the Campbell-Hausdorff series
corresponding to the Lie structure constants $\gm^\mu_{\rho\si}$ in  the
Lie algebra ${\bf V}^*$, and ${\al \cdot X}$ for the
matrix with entries $\al^i_{\rho k} X^\rho$.
The coproduct is evidently coassociative,
as the elements $X^\mu$ commute. The problem boils down  to
evaluating the full set of quantum commutation relations
consistent with (\ref{coproduct:1}).
We will search for them in the form
\begin{equation}
\begin{array}{lcl}
 [X^\mu,X^\nu] &=& 0,\\[6pt]
 [H_i, H_k]&=& C(X)^m_{ik}H_m , \\[6pt]
 [H_i,X^\mu] &=& A(X)_i^\mu,\\[6pt]
\end{array}\label{COMMREL1}
\end{equation}
treating quantum structure constants as formal series in
commutative generators $X^\mu$.
\begin{theor}
There exists the unique
quantization of the bialgebra $({\bf L},{\bf L}^*)$
with coproduct (\ref{coproduct:1})
and commutation relations (\ref{COMMREL1}), such
that
$$
C(0)^m_{ik}=C^m_{ik}, \quad
\frac{\partial A(0)_i^\mu}{\partial X^\nu}=A_{i\nu}^{\>\>\mu}.
$$
It is a functor from the category ${\cal B}$
onto the sub-category ${\cal H}$ of Hopf algebras.
\label{main}
\end{theor}
\underline{Proof.}
Substituting (\ref{coproduct:1}) into $[\Delta(H),\Delta(X)]=\Delta([H,X])$
we come to the equation
\be
 A(D(X',X''))_i^\mu &=&
 (e^{\al \cdot X'})^k_i
 \partial''_\nu D^\mu(X',X'') A(X'')_k^\nu+
 \partial'_\nu D^\mu(X',X'') A(X')_i^\nu,
 \label{eq:1}
\ee
where primes distinguish tensor factors.
Regarding $X^\mu$ as the coordinate functions on
the Lie group $\exp({\bf V}^*)$ we can consider
$A(X)_i^\mu$ as a set of vector fields
in the normal neighborhood of the identity,
labeled by index $i$.
Then equation (\ref{eq:1}) is nothing else
than
\be
 A(\xi\circ \z)_i^\mu &=&
 (e^{\al \cdot  \xi})^k_i  L_\xi A(\z)_k^\mu+
 R_\z A(\xi)_i^\mu, \quad \xi,\z \in {\bf V}^*,
 \label{eq:2}
\ee
where $L_\xi$, $R_\xi$ stand for the  left and right
actions of the group $\exp({\bf V}^*)$ on the vector
fields. Note that both these actions
commute with the action specified by the matrices
$\al_\mu$. Transition  to the functions $\hat{A}(\xi)=R^{-1}_\xi A(\xi)$
leads to the group 1-cocycle equation
\be
 \hat{A}(\xi\circ \z)_i^\mu &=&
 (e^{\al \cdot  \xi})^k_i  Ad(\xi) \hat{A}(\z)_k^\mu+
  \hat{A}(\xi)_i^\mu,
 \label{eq:3}
\ee
which has the unique solution, provided the differential $d\hat{A}(0)$
is a corresponding 1-cocycle of the Lie algebra ${\bf V}^*$.
That is a part of the Lie bialgebra consistency conditions (\ref{cond1}) on
the pair $({\bf L},{\bf L}^*)$. The explicit formula for
the functions $A(X)_i^\nu$ is
\be
A(X)_i^\mu=\biggl(\frac{\scriptstyle \gm''\cdot X}{e^{\gm''\cdot X}-1}\quad
\frac{ e^{\al'\cdot X + \gm''\cdot X}-1}{\scriptstyle \al'\cdot X +
\gm''\cdot X}\biggr)^{k\mu}_{i\nu}
A_{k\rho}^{\>\>\nu}X^\rho .
\label{expl}
\ee
Here $(\gm\cdot X)^\mu_\nu = \gm^\mu_{\si\nu} X^\si$ specifies the adjoint
representation of the Lie algebra ${\bf V}^*$.
We mark the matrices with primes to stress that they act on the different
groups of indices.
Note, that  formula (\ref{expl}) is simplified in the case of
Abelian ${\bf V}^*$: then $A(X)_i^\mu$ takes the form
$
A(X)_i^\mu=\biggl(\frac{e^{\al\cdot X} -1}{\scriptstyle \al\cdot X}\biggr)^\mu_\nu
A_{i\rho}^{\>\>\nu}X^\rho.
$

Requirement $[\Delta(H),\Delta(H)]=\Delta([H,H])$
leads to the following two equations:
\be
C(D(X',X''))^k_{ij}=C(X')^k_{ij},
\ee
meaning that $C(X)^i_{jk}$ are actually constant, and
\be
(e^{\al \cdot  X})^k_mC^m_{ij}-
C^k_{mn}(e^{\al \cdot  X})^m_i(e^{\al \cdot  X})^n_j
&=&
[H_i,(e^{\al \cdot  X})^k_j]-[H_j,(e^{\al \cdot  X})^k_i].
 \label{eq:4}
\ee
Lie algebra representation by matrices $\al_\mu$ induces
an anti-homomorphism of $\exp({\bf V}^*)$ into the
linear group {\cal Lin}$({\bf H},{\bf H})$.
The expressions on the right-hand side of (\ref{eq:4}) are
the vector fields $A(X)_i$ transferred by that map to
{\cal Lin}$({\bf H},{\bf H})$. In terms of matrices
$a = e^{\al \cdot X}$, we rewrite (\ref{eq:4}) as
\be
a^k_m C^m_{ij}-
C^k_{mn}a^m_i a^n_j
&=&
A(a)_{ij}^{\>\>k}-
A(a)_{ji}^{\>\>k}
 \label{eq:5}
\ee
or
\be
a C (a^{-1}\oo a^{-1}) - C
&=&
A(a)\pi (a^{-1}\oo a^{-1}),
 \label{eq:6}
\ee
where we introduced the anti-symmetrizer
$\pi$, $\pi^{ij}_{kl}=\delta^i_k\delta^j_l - \delta^i_l\delta^j_k$.
The left-hand side of this equation is a coboundary
1-cocycle on the linear group, so we must prove that for the
right-hand side. Then, since group 1-cocycles are uniquely
determined by their derivatives at the identity, (\ref{eq:6}) will
follow from  (\ref{cond2}).
By virtue of (\ref{eq:2}), we have
\be
A(ba)\pi ((ba)^{-1}\oo (ba)^{-1}) &=&
 (A(b)(a\oo a) + b A(a))(a^{-1}\oo a^{-1})\pi (b^{-1}\oo b^{-1})
 \n\\&=&
 A(b)\pi (b^{-1}\oo b^{-1}) + b \{A(a)\pi (a^{-1}\oo a^{-1})\} (b^{-1}\oo
b^{-1}),
\n
\ee
as required.

The counit is evident: $\epsilon(H_i)=\epsilon(X^\mu) = 0$.
The antipode is determined on the generators by the coproduct:
$S(X^\mu)=-X^\mu$,
$S(H_i)=-(e^{-\al\cdot X})^k_i H_k$. Let us prove that it is
extended over the whole algebra anti-homomorphically. It is
trivial in the commutative $X$-sector. Further,
\be
[S(H_i),S(X^\mu)]&=&(e^{-\al\cdot X})^k_i[ H_k,X^\mu]=
(e^{-\al\cdot X})^k_i A(X)_k^\mu
\n\\&=&
-A(-X)_i^\mu=S([X^\mu, H_i]),\n
\ee
as  immediately follows from formula (\ref{expl}).
Condition $[S(H_i),S(H_j)]= S([H_j,H_i])$ boils down
to verification of
$$
C^k_{mn}{a^{-1}}^{m}_{i} {a^{-1}}^{n}_{j}
+{a^{-1}}^{m}_{i}[H_{m},{a^{-1}}^{n}_{j}]
+{a^{-1}}^{n}_{j}[{a^{-1}}^{m}_{i},H_{n}]=
{a^{-1}}^k_m C^m_{ij}, \quad a=e^{\al\cdot X}.
$$
We represent it as
$$
S\Bigl(C^k_{mn}a^{m}_{i} a^{n}_{j}\Bigr)
-S\Bigl([H_{i},a^{n}_{j}]\Bigr)
-S\Bigl([a^{m}_{i},H_{j}]\Bigr)=
S\Bigl(a^k_m C^m_{ij}\Bigr),
$$
which holds true in view of (\ref{eq:5}).

Thus we described the Hopf structure of quantum algebras from
${\cal H}$.
We have  yet to check that the quantization ${\cal B} \to {\cal H}$
is a natural map of categories. Let $\phi$ be a Lie bialgebra morphism
${\bf L}\to {\bf L}'$ such that  $\phi\si= \si'\phi$.
This implies that $\phi({\bf H})\subset{\bf H}'$ and
$\phi({\bf V})\subset{\bf V}'$.
We define the map $\Phi\colon U_q({\bf L})\to U_q({\bf L}')$
by the same formulas on the generators as $\phi$. Note that
$X^\mu$ and $H_i$ in the quantum algebra are, normally, not
the same as the classical generators. We do not interested in relations
between them, although in some cases like twisted algebras
\cite{M12,KLM} it is possible to give the explicit formulas.
Linear map $\Phi$ on the generators can be extended over the whole
quantum algebras as a Hopf homomorphism. It is evident for
the coproduct because it is given by the composition in
the dual Lie groups, and $\phi$ is a Lie bialgebra homomorphism.
That can be shown for the commutation relations as well.
Indeed, value of the quantum commutator (\ref{expl}) differs
from the classical one by involvement of the matrices
$(\al\cdot X)^i_k$ and $(\gm\cdot X)^\mu_\nu$.
They specify the adjoint representation of the subalgebra
${\bf V}^*\subset {\bf L}^*$. Because $\phi^*$ is a homomorphism
of the dual Lie algebras, and preserves $\si$-invariant subspaces,
matrices $\Phi^i_k$ and $\Phi^\mu_\nu$ are pulled through
$\al\cdot X$ and $\gm\cdot X$ properly, e.g.
$ \Bigl(\al'\cdot X'\Bigr) \Phi = \Phi\Bigl(\al\cdot\Phi(X)\Bigr) $,
so the proof becomes immediate.

We denote quantization of $U({\bf L})$ as $U_q({\bf L})$ although
there is no deformation parameter involved so far. It can be introduced
by substitution $\al \to \ln(q) \al$, $\gm \to \ln(q) \gm$ for the structure
constants of the dual Lie algebra but irrelevant for
our study. Algebra $U_q({\bf L})$ contains two classical objects:
universal enveloping algebra $U({\bf H})$ and the commutative
algebra of functions on the Lie group $\exp({\bf V}^*)$.
In accordance with our convention, we may
assume $\mbox{Fun}(\exp({\bf V}^*))\sim U_q({\bf V})$.
Actually
$U_q({\bf L})$ is a bicrossproduct Hopf algebra
$U({\bf H})\triangleright U_q({\bf V})$, with
the coaction on $U({\bf H})$ given by
$H_i\to (e^{\al \cdot X})^k_i \oo H_k $.

\section{Duality- and double-invariance}
\label{DDI}
Lie bialgebra category ${\cal B}$ is evidently self-dual, the involution
$\si$ going
over into $-\si^*$. Let us prove the analogous assertion for
${\cal H}$ and  deduce explicitly the canonical element.

As a linear space $U_q({\bf L})$ is the tensor product
$U_q({\bf V})\oo U({\bf H})$. There are two natural algebra maps
from $U^*({\bf H})$ and $U({\bf V^*})$
into $U^*_q({\bf L})$: we set $\h\to\ve_{\bf V}\oo\h$ and
$\z\to\z\oo\ve_{\bf H}$, correspondingly. It is straightforward that
$$ \langle \h \z,\varphi(X)\psi(H)\rangle =
   \langle \h\oo \z,\varphi(\Delta(X))
   \psi(\Delta(H))\rangle =
   \langle \h ,\psi(H)\rangle  \langle \z,\varphi(X)\rangle.
$$
Here we do not make difference between functionals $\z$, $\h$ and
there images in $U^*_q({\bf L})$. Just proved, the  factorization
property justifies such an abuse of notations. It means that
linear spaces $U^*({\bf H})$ and $U({\bf V^*})$ are isomorphically
embedded into $U^*_q({\bf L})$ (in fact, these are isomorphisms of
associative algebras, see Appendix), and the induced
map $U({\bf V^*})\oo U^*({\bf H}) \to U({\bf V^*})U^*({\bf H})$
is a linear bijection on $U^*_q({\bf L})$.

Let us choose the bases $\h^i\in U^*({\bf H})$
and $\z_\mu \in U({\bf V^*})$ dual to $H_i$ and $X^\mu$ as generators
of $U^*_q({\bf L})$.
It can be shown that they have coproducts of the form (\ref{coproduct:1}),
and commutation relations
similar to  (\ref{COMMREL}), of course,  after interchanging ${\bf L}$ and ${\bf L^*}$.
That is done in Appendix.
Then Theorem \ref{main} states that  ${\bf L}^*$ admits the unique
quantization $U_q({\bf L}^*)\sim U^*_q({\bf L})$ belonging to
${\cal H}$.

Because of the factorization property and due to the fact that
the pairings $\langle \h ,\psi(H)\rangle $
and $\langle \z,\varphi(X)\rangle$ are the same as if
 $H_i$ and  $\z_\mu$ were primitive, we can easily
 write down the canonical element
 ${\cal T}\in  U_q({\bf L}^*)\oo U_q({\bf L})$.
Nevertheless, it is convenient to deal with the
opposite algebra $U_q({\bf L}^*)_{op}$; moreover, it
is that algebra which takes part in construction of the double,
the subject of our  further interest.
In Appendix we prove the following result:
\be
{\cal T}= \exp(\z_\mu \oo X^\mu )\exp(\h^i\oo H_i),
\label{T-matrix}
\ee
using expressions for the canonical elements
of the classical universal enveloping Hopf algebras. Summation over
repeating indices is assumed.

Now we proceed to the study of the double in the category ${\cal H}$.
First recall that the double of a Lie bialgebra ${\cal D}({\bf L})$
is a unique Lie bialgebra such that ${\bf L}$ and ${\bf L}^*$, taken
with the opposite Lie bracket, are embedded as
sub-bialgebras, and the canonical pairing between them gives rise to
a non-degenerate invariant symmetric bilinear form on
${\cal D}({\bf L})$. The double procedure preserves category ${\cal B}$.
Indeed, the classical commutation relations followed from the
definition are
\begin{equation}
\begin{array}{lcl}
 [H_i, H_k]&=& C^m_{ik}H_m ,
 \\[6pt]
 [\z_\mu,\z_\nu]   &=&\gm^\si_{\nu\mu}\z_\si,
 \\[6pt]
 [H_i,\z_\mu]   &=&-\al_{\mu i}^{\>\>k}H_k-A_{i\mu}^{\>\>\nu}\z_\nu
 \\[6pt]
 [X^\mu,X^\nu]     &=&0,
 \\[6pt]
 [\h^i,\h^j]       &=&0,
 \\[6pt]
 [X^\mu,\h^i]       &=&0,
 \\[6pt]
 [H_i,X^\mu] &=& A_{i\nu}^{\>\>\mu} X^\nu,
 \\[6pt]
 [H_i,\h^j] &=&  C^j_{ki} \h^k +  \al_{\mu i}^{\>\>j} X^\mu,
 \\[6pt]
 [\z_\mu,\h_i]   &=&-\al_{\mu i}^{\>\>k}\h_k,
 \\[6pt]
 [\z_\mu,X^\nu]  &=& -\gm^\nu_{\si\mu}X^\si - A_{i\mu}^{\>\>\nu}\h^i,
\end{array}\label{COMMRELdouble}
\end{equation}
that proves the assertion.

It is thus natural to expect the analogous statement
in the quantum case.
\begin{theor}
Quantum double preserves category ${\cal H}$.
Moreover,
${\cal D}(U_q({\bf L}))=U_q({\cal D}({\bf L}))$.
\label{thlast}
\end{theor}
As a coalgebra, the double coincides with the tensor
product of $U_q({\bf L})$ and $U_q({\bf L}^*)_{op}$,
which are at the same time subalgebras.
Therefore, to prove the theorem, it suffices to show that
the cross-relations have the appropriate form.
Then we will satisfy conditions of Theorem \ref{main}
which states the uniqueness of the quantization and
provides its explicit form.
The cross-relations are deduced from the Yang-Baxter equation
on the canonical element and can be written
as
$$
e_\mu e^\nu  =e^\al e_\bt m^{\>\>\nu}_{\gm\al\si}
m_{\>\>\mu}^{\rho\bt\si}S_\rho^\gm,
$$
where $e_\mu \in U_q({\bf L})$, $e^\mu \in U_q({\bf L}^*)$,
$m^{\>\>\nu}_{\gm\al\si}$ and $ m_{\>\>\mu}^{\rho\bt\si}$ denote
the  iterated
coproduct structure constants, and $S_\rho^\gm$ is the matrix of
the antipode.
Using the explicit formulas for the coproducts
\begin{equation}
\begin{array}{lcl}
 (\Delta\oo id)\circ\Delta (X) &=& D(X\oo 1\oo 1,1\oo X\oo 1, 1\oo1\oo X),\\
 (\Delta\oo id)\circ\Delta (H) &=& e^{\al \cdot X} \oo e^{\al \cdot X} \oo H +
                  e^{\al \cdot X} \oo H \oo 1 + H  \oo 1 \oo 1,\\
 (\Delta\oo id)\circ\Delta (\h) &=& D(\h\oo 1\oo 1,1\oo \h\oo 1, 1\oo1\oo \h),\\
 (\Delta\oo id)\circ\Delta (\z) &=& e^{A \cdot \h} \oo e^{A \cdot \h} \oo \z   +
                  e^{A \cdot \h} \oo \z \oo 1 + \z  \oo 1 \oo 1,
\end{array}
\end{equation}
we  get the required result. For example,
$$
H\z = \langle e^{-\al \cdot X},\z \rangle H+
\langle e^{-\al \cdot X},e^{A \cdot \h} \rangle \z H +
\langle H,e^{-A \cdot \h} \rangle \z
$$
(only non-vanishing terms retained). Thus we
obtain the expression for the commutator
$[H_i,\z_\mu]  =-\al_{\mu i}^{\>\>k}H_k-A_{i\mu}^{\>\>\nu}\z_\nu.$
Similarly we prove that $X^\mu$ and $\h^i$ form a
commutative algebra, invariant under the  adjoint action of
$H_i$  and $\z_\mu$.

Now consider a Lie bialgebra from ${\cal B}$
which possesses a classical r-matrix. To use  the advantages
of functorial property of the quantization, let us
demand that the r-matrix viewed as a Lie bialgebra morphism
be that in ${\cal B}$. It means that $r$ must
commute with the involution $\si$ (we
remind that for dual bialgebra $\si$ goes over into
$-\si^*$). The general form of $r$ is then
$$r={\cal P}^i_\mu H_i\oo X^\mu + {\cal Q}_\mu ^i X^\mu \oo H_i.$$
We may assume that summation is performed over
elements of the bases of  ${\bf L}_+$ and ${\bf L}_-$.
Then the universal R-matrix is found from
(\ref{T-matrix}):
$$
{\cal R}=\exp({\cal P}^i_\mu H_i\oo X^\mu  )
         \exp({\cal Q}_\mu ^i X^\mu\oo H_i).
$$
In this form this is a generalization of the result of \cite{BH}
obtained for the simplest case of the double of the jordanian
Borel quantum algebra.

\section{Discussions}
Drinfeld's conjecture of the possibility to quantize an arbitrary Lie
bialgebra was
proved by Etingof and Kazhdan \cite{EK}. Although there are numerous
examples of quantum algebras, the problem of exponentiating
a Lie bialgebra in  every particular case remains highly nontrivial.
In the present paper, we do it for a class of
algebras playing significant role in the
classical differential geometry and physics, inhomogeneous
Lie algebras. This class forms a nice
category invariant under dualization and quantum
double operations. The quantization built is a functor from
the Lie bialgebra category of concern into the category of Hopf
algebras. This functor commutes with the functors of dualization
and double. We showed that the quantization of inhomogeneous
Lie algebras possessing classical r-matrix itself contains the
solution to the quantum Yang-Baxter equation which
is the universal R-matrix for its quasitriangular
Hopf subalgebra. This statement is based on the fact that
 every minimal quasitriangular Hopf algebra is a quotient
of the double of its subalgebra, and, also, on the functorial property
of the quantization.

The class of Hopf algebras considered in the present article
includes twisted universal enveloping Lie algebras
taking part in the null-plane quantization of Poincar\'e
algebra and extended jordanian deformations of $sl(N)$.
(\cite{M1,M12,KLM}). They are characterized by the identification
${\bf V}\sim {\bf H}^*$ and involve non-degenerate
Lie algebra 1-cocycles in building crossproduct
${\bf H}\triangleright{\bf H}^*$. The resulting
quantum algebras were triangular and twist-equivalent
to the classical universal enveloping algebras.
Here we have studied the general case.

The appearance of Lie group 1-cocycles in construction of quantization
is quite understandable. According to Drinfeld \cite{D3}, a group
1-cocycles with values in Lie algebra external square
defines  Poisson-Lie structures on the group.   Generators
$H_i \in {\bf H}$ and  $X^\mu \in {\bf V}$ are the coordinate
functions on the group ${\bf H}^* \triangleleft\exp({\bf V}^*)$,
and the quantum commutation relations among them represent nothing else than
the Poisson bracket. Along with the coproduct (\ref{coproduct:1}),
this implies that the quantum algebra $U_q({\bf L})$ is just
universal enveloping algebra of the Lie-Poisson algebra on the
group $\exp({\bf L}^*)$. Indeed, the classical commutation relations
of the types $[H,H]$ and  $[X,X]$ are given by the Poisson structure
on the Abelian group ${\bf H}^*$ and the trivial Poisson bracket
on the group $\exp({\bf V}^*)$. Further, the Poisson bracket of the
type $\{H,X\}$ must satisfy the relation
\begin{equation}
\begin{array}{lcl}
\Delta(\{H,X\})&=& \{e^{\al\cdot X}\oo H +H\oo 1,D(X\oo 1,1\oo X) \}
\\&=&
(e^{\al\cdot X}\oo 1) \{1\oo H,D(X\oo 1,1\oo X) \}
\\&+&
\{H\oo 1,D(X\oo 1,1\oo X) \}
\\
\Delta(\{H,H\})&=&
\{e^{\al\cdot X}\oo H +H\oo 1,
e^{\al\cdot X}\oo H + H\oo 1 \}=
\\&=&
e^{\al\cdot X} e^{\al\cdot X}\oo\{ H,H \}+
\{ H, H \}\oo 1+
\\&+&
\{e^{\al\cdot X},H\}\oo H +
\{H,e^{\al\cdot X} \}\oo H
\end{array}
\label{poisson}
\end{equation}
(we drop all the indices for the reason of transparency).
These expressions involve the Poisson bracket
and the multiplication in the function algebra
on $\exp({\bf L}^*)$. It is seen that  only the
product in
Fun($\exp({\bf V}^*)$) really matters, which is commutative
and coincides with that on $U_q({\bf V})\subset U_q({\bf L})$. So,
one can consider equations (\ref {poisson})
as in the quantum algebra.

\vspace{1cm}
\noindent
{\Large \bf Acknowledgement}

\vspace{0.5cm}
\noindent
We are grateful to T. Hodges, N. Yu. Reshetikhin, and
M. A. Semenov-Tyan-Shansky for illuminating discussions and
helpful remarks.

\vspace{1cm}
\noindent
{\Large \bf Appendix}

\vspace{0.5cm}
\noindent
The aim of this section is to exhibit the details of the proof
of the formula (\ref{T-matrix}).
The canonical element for the universal enveloping algebra
$U({\bf V}^*)$ is
\cite{D3}
$$e^{\z_\mu \oo X^\mu}= \sum_n
\sum_{\vec \mu}(\z_{\mu_1},\ldots,\z_{\mu_n})\oo X^{\mu_1}\ldots X^{\mu_n},
$$
where $\vec \mu =({\mu_1},\ldots,{\mu_n})$ stands for the ordered
multiindex of length $n$, and
the parentheses denote symmetrized monomials
$
(\z_{\mu_1},\ldots,\z_{\mu_n})={1\over n !}{1\over s(\mu)}\sum_{\si}
 \z_{\si(\mu_1)} \ldots \z_{\si(\mu_n)}
$
of degree $n$. Here  $\si$ belongs to the symmeric group $S_n$
and $s(\mu)$ is equal to the order of the stability subgroup
under permutations of $\mu$.
Similarly we have
$$
e^{\h^i \oo H_i}= \sum_n
\sum_{\vec i}\h^{i_1}\ldots\h^{i_n}\oo (H_{i_1},\ldots ,H_{i_n})
$$
for the algebra  $U({\bf H})$.
Hence, due to the factorization of the matrix elements of the canonical
pairing (Section \ref{DDI}), the  element
$e^{\z_\mu \oo X^\mu}e^{\h^i \oo H_i}$ is canonical for
$U_q({\bf L})^*_{op}$. What remains is to show that the
coproduct and commutation relations for
$\z_\mu$  and for $\h^i$ have the proper form
\be
\langle  \h^i \h^j, X^{\mu_1}\ldots X^{\mu_n}(H_{i_1},\ldots ,H_{i_m})\rangle &=&
\ve(X^{\mu_1}\ldots X^{\mu_n})
\langle \h^i\oo \h^j,\Delta((H_{i_1},\ldots ,H_{i_m}))\rangle \n\\&=&
\ve(X^{\mu_1}\ldots X^{\mu_n})
\langle \h^i\oo \h^j,\Delta_0((H_{i_1},\ldots ,H_{i_m}))\rangle .\n
\ee
The last equality, where $\Delta_0$ is the classical cocommutative
comultiplication in  $U({\bf H})$,
is due to that one can push all the factors
$(e^{\al \cdot X})^k_i$ in  $\Delta(H_{i_k})$ to the left, where
they are reduced to $1$,  as
though they commute
with $H$'s. That is because $X$'s generate an ideal, and, once appeared,
the terms containing them will be annihilated by $\h^i$.
So, the  generators $\h^i$ commute.
Further,
\be
\langle  \z_\mu \z_\nu, X^{\mu_1}\ldots X^{\mu_n}(H_{i_1},\ldots ,H_{i_m})\rangle &=&
\langle \z_\mu \oo\z_\nu,\Delta(X^{\mu_1}\ldots X^{\mu_n}(H_{i_1},\ldots ,H_{i_m}))  \rangle \n\\&=&
\langle \z_\mu \oo\z_\nu,\Delta(X^{\mu_1}\ldots X^{\mu_n})((H_{i_1},\ldots ,H_{i_m})\oo 1)  \rangle \n\\&=&
\ve(H_{i_1},\ldots ,H_{i_m})
\langle \z_\mu \oo\z_\nu,\Delta(X^{\mu_1}\ldots X^{\mu_n}) \rangle, \n
\ee
and therefore $\z$'s have the classical product of the universal enveloping
algebra $U({\bf V}^*)$.
Among the matrix elements
\be
\langle  \z_\mu \h^i, X^{\mu_1}\ldots X^{\mu_n}(H_{i_1},\ldots ,H_{i_m})\rangle &=&
\langle \z_\mu \oo \h^i,(X^{\mu_1}\ldots X^{\mu_n}\oo1)\Delta(H_{i_1},\ldots ,H_{i_m})  \rangle \n
\ee
only those survive where $n\leq 1$. Developing products of $\Delta(H)$'s
we see that monomials $1\oo H_{i_1}\ldots H_{i_k}$ turn out to be symmetrized automatically,
hence we can retain terms of the first degree in $1\oo H_i$ only.
And furthermore, if $n=1$ then with neccessaty $m=1$.
The non-vanishing pairings are
\be
\langle  \z_\mu \h^i, X^{\mu_1}(H_{i_1},\ldots ,H_{i_m})\rangle &=&
\langle \z_\mu \oo \h^i,X^{\mu_1}(H_{i_1},\ldots ,H_{i_{m-1}})(e^{\al \cdot X})^k_{i_m} \oo H_k   \rangle
\n\\
 &=&
\langle \z_\mu \oo \h^i,\varphi(X) \oo H_k   \rangle \n
\ee
where $\varphi(X)$ is a result of pulling the exponential to the left.
Thus we state that the difference $\z_\mu \h^i - \h^i\z_\mu$
does not vanish only on the elements $(H_{i_1},\ldots ,H_{i_m})$, therefore
it depends solely on $\h^i$.

We have yet to find the coproduct. It is straightforward
that $\Delta(\h)$ survives on the elements with no $X$'s and
therefore is expressed by the Campbell-Hausdorff series
corresponding to $U({\bf H})$. For $\Delta(\z)$ the nontrivial pairing
is with elements $  X^\mu\oo 1$ and $(H_{i_1},\ldots ,H_{i_m})\oo X^\mu$. While
pulling $H$'s to the right we can assume that they commute
with $X$'s with the classical relations, that is the
commutator is linear in $X$ because the higher degrees will annihilate.
Thus we come to the desired formula
$$ \Delta(\z_\mu)=(e^{A \cdot \h})^\nu_\mu \oo \z_\nu + \z_\mu \oo 1.$$


\begin{thebibliography}{18}
\bibitem{Mj4} S. Majid,
{\em Crossproduct quantization, nonabelian cohomology and twisting of Hopf
    algebras},
    hep-th/9311184.
\bibitem{CK} J. A. de Azc\'arraga, M. A. del Olmo, J. C. P\'erez Bueno,
and M. Santander,
{\em Graded contractions and bicrossproduct structure of deformed
inhomogeneous algebras}, J. Phys. A {\bf 30} (1997) 3069
(q-alg/9612022);

S. Majid and H. Ruegg,
{\em Bicrossproduct structure of $\kappa$-Poincar\'e group and
non-commutative  geometry},
   Phys. Lett. B {\bf 334} (1994) 348 (hep-th/9405107).
\bibitem{Ogiv} O. Ogievetski,
{\em Hopf structures on the Borel subalgebra of $sl(2)$},
 Max-Plank-Institute preprint VPI--Ph/92--99.
\bibitem{D1} V.G. Drinfeld,
{\em On constant quasiclassical solutions to the quantum Yang-Baxter
equation},
   DAN USSR  {\bf 273} \# 3 (1983) 531.
\bibitem{D2} V.G. Drinfeld,
{\em Quasi-Hopf algebras},
   Leningrad Math. J. {\bf 1} (1990) 1419.
\bibitem{AHO} O. Arratia, F. J. Herranz, and M. A.  del Olmo,
{\em  Bicrossproduct structure of the null-plane quantum  Poincar\'e
algebra},
 J. Phys. A {\bf 31} 1998  L1.
\bibitem{M12} A. I. Mudrov,
    {\em Twisting cocycles in fundamental representation and
    triangular bicrossproduct Hopf algebras},
    J. Math. Phys. {\bf 39}  (1998) 5608 (math.QA/9804024).
\bibitem{KLM} P. P. Kulish, V. D. Lyakhovsky, A. I. Mudrov,
    {\em Extended jordanian twists for Lie algebras},
    J. Math. Phys., to appear (math.QA/9806014).
\bibitem{ES} P. Etingof, A. Soloviev,
    {\em Quantization of geometrical classical r-matrices}\\
    (math.QA/9811001).
\bibitem{Mj2} E. Beggs, S. Majid,
{\em Quasitriangular and differential structures on bicrossproduct
    Hopf algebras},
    q-alg/9701041.
\bibitem{Rad} D. E. Radford
{\em Minimal quasitriangular Hopf algebras},
    J.Algebra {\bf 157} (1993) 285.
\bibitem{D3} V. G. Drinfeld,
{\em Quantum Groups}, in Proc. Int. Congress of
Mathematicians, Berkeley, 1986, ed. A.V. Gleason, AMS, Providence, (1987) 798.

\bibitem{RS}  N. Yu. Reshetikhin, M. A. Semenov-Tian-Shansky,
{\em Quantum $R$-matrices and factorization problem},
  J. Geom. Phys. {\bf 5} (1988) 533.

\bibitem{BH} C. Burdik, P. Hellinger,
           {\em The quantum double for a nonstandard
           deformation of a Borel subalgebra $sl(2,C)$}
           (hep-th/9303035).
\bibitem{V} A. A. Vladimirov,
             {\em On quasitriangular Hopf algebras related to the
              Borel subalgebra of $sl_2$},
              Symm.Meth.in Phys., {\bf 2}, Dubna (1994) 574.
\bibitem{STS9} M. A. Semenov--Tian--Shanski,
         {\em  Poisson--Lie groups,
        the quantum duality principle and the twisted quantum duble},
        Theor. Math. Phys.  {\bf 93}  \# 7 (1992) 302.
\bibitem{LM} V. D. Lyakhovsky,
    {\em The role of solvable groups in quantization of Lie algebras}
    Zapiski nauch. sem. POMI {\bf 209} (1993) 131 (hep-th/9405044);
     A. I. Mudrov,
    {\em Quantum Lie algebras with dual groups of the simplest type},
    Vest.SPbSU {\bf 4} \# 4 (1994)  3.
\bibitem{EK}  P. Etingof, D. Kazhdan,
    {\em Quantization of Lie bialgebras},
    Selecta Math. {\bf 2}, \# 1  (1996) 1 (q-alg/9510020).

\bibitem{M1} A. I. Mudrov,
{\em Twisting cocycle for null-plane quantized Poincar\'e algebra},
   J. Phys. A {\bf 31} (1998), 6219 (q-alg/9711001).

\end{thebibliography}
\end{document}